# SUBGEOMETRIC ERGODICITY OF STRONG MARKOV PROCESSES

By G. Fort and G. O. Roberts

*CNRS/LMC-IMAG and Lancaster University*

We derive sufficient conditions for subgeometric $f$-ergodicity of strongly Markovian processes. We first propose a criterion based on modulated moment of some delayed return-time to a petite set. We then formulate a criterion for polynomial $f$-ergodicity in terms of a drift condition on the generator. Applications to specific processes are considered, including Langevin tempered diffusions on $\mathbb{R}^n$ and storage models.

**1. Introduction.** This paper is devoted to the study of subgeometric $f$-ergodicity of a strong Markov semigroup $(P^t)_{t\geq 0}$. That is, for a subgeometrically increasing rate function $r := (r(t))_{t\geq 0}$, and a Borel function $f \geq 1$, we propose sufficient conditions implying the limit

$$\lim_{t\to+\infty} r(t)\|P^t(x,\cdot) - \pi(\cdot)\|_f = 0,$$

for $\pi$-almost all (a.a.) $x$, where $\pi$ is the unique invariant probability measure. Our main condition is couched in terms of modulated moments of return-times to a "test-set." In this form, this condition extends earlier criteria, implying different notions of stability (such as Harris-recurrence, positive Harris-recurrence, ergodicity and $f$-ergodicity) for continuous-time Markov processes. This condition is also analogous to the criterion for subgeometric $f$-ergodicity of discrete time Markov chains. We also derive a condition for polynomial ergodicity which is easy to check in many applications. This condition is expressed in terms of inequality on the semigroup *generator*, and is analogous to the so-called *drift inequality* in the discrete-time case.

We apply our results to the study of strongly Markovian processes, giving three nontrivial examples, two of which are of considerable applied probabilistic interest. We first consider a simple jump process as a toy example,

---

demonstrating that $f$-ergodicity at a logarithmic (resp. polynomial or subexponential) rate is narrowly related to the existence of a logarithmic (resp. polynomial or sub-exponential) moment of the mean-time spent in each state, with respect to the jump distribution. We then consider Langevin tempered diffusions on $\mathbb{R}^n$ which are relevant to Markov Chain Monte Carlo (MCMC) techniques since they construct a diffusion process with given stationary distribution $\pi$ (which only needs to be available up to an unknown normalization constant). When the stationary distribution is polynomial in the tails, the (simple) Langevin diffusion can not be ergodic at a geometric rate and we show that it is polynomially ergodic. We also consider Langevin tempered diffusion in which the diffusion matrix is a scalar matrix with coefficient $\pi^{-2d}$, $d > 0$, and prove that even when the target distribution is polynomial in the tails, a convenient choice of the temperature $d$ involves geometric ergodicity of the process. Finally, we study a compound Poisson-process driven Ornstein–Uhlenbeck process which is used in storage models and more recently in financial econometrics. It is known that when the distribution of the jump $F$ has sufficiently light tails, the process is geometrically ergodic. We investigate the case where $F$ is heavy tailed and establish the sub-geometric ergodicity of the process under appropriate conditions in this case.

The paper is organized as follows. We first recall basic definitions on Markov process, as well as reviewing existing results on ergodicity of strongly Markovian process. The new criteria for subgeometric ergodicity are given in Section 2, and the proofs are postponed to Section 4. Section 3 is devoted to the three examples mentioned above.

### 1.1. *Basic definitions on Markov process.*

Let $\mathcal{X}$ be a locally compact and separable metric space endowed with the Borel $\sigma$-field $\mathcal{B}(\mathcal{X})$. $X = (\Omega, A, (\mathcal{F}_t)_{t \geq 0}, (X_t)_{t \geq 0}, \mathbb{P}_x)$ is a $\mathcal{X}$-valued Borel right process so that it is a temporally homogeneous Markov process, strongly Markovian with right-continuous sample paths (see, e.g., [31]). $\mathbb{P}_x$ (resp. $\mathbb{E}_x$) denotes the canonical probability (resp. expectation) associated to the Markov process with initial distribution $\delta_x$, the Dirac distribution at point $x$. Let $(P^t)_{t \geq 0}$ be the associated Markov semigroup.

We recall basic definitions and properties on Markov process that will be used throughout this paper. The process $X$ is *$\phi$-irreducible* for some $\sigma$-finite measure $\phi$ on $\mathcal{B}(\mathcal{X})$ if

$$\phi(A) > 0 \quad \Longrightarrow \quad \mathbb{E}_x\left[\int_0^\infty \mathbb{1}_A(X_s)\,ds\right] > 0 \qquad \forall x \in \mathcal{X}.$$

If the process is $\phi$-irreducible, there exists a maximal irreducibility measure $\psi$ that dominates any irreducibility measure [24]. In fact, if $\pi$ is an invariant measure, that is, $\pi P^t = P^t$ for all $t \geq 0$, then $\pi$ is a maximal irreducibility measure. Any measurable set which is of positive $\psi$-measure is said to be *accessible*. A set $C \in \mathcal{B}(\mathcal{X})$ is $\nu_b$-petite for the process (or simply petite) if there exists a probability measure $b$ (resp. nontrivial $\sigma$-finite measure $\nu_b$) on the Borel $\sigma$-field of $\mathbb{R}_+$ [resp. on $\mathcal{B}(\mathcal{X})$] such that

$$\int_0^\infty P^t(x, \cdot)b(dt) \geq \nu_b(\cdot) \qquad \text{for all } x \in C.$$

A $\phi$-irreducible process always possesses an accessible closed petite set ([20], Proposition 3.2). A process is *Harris-recurrent* if there exists a $\sigma$-finite measure $\phi$ such that

$$\phi(A) > 0 \quad \Longrightarrow \quad \mathbb{P}_x\left(\int_0^\infty \mathbb{1}_A(X_s)\,ds = +\infty\right) = 1, \qquad x \in \mathcal{X};$$

or, equivalently, if there exists a $\sigma$-finite measure $\mu$ such that $\mu(A) > 0 \Longrightarrow \mathbb{P}_x(\tau_A < \infty) = 1$ for all $x \in \mathcal{X}$ where $\tau_A$ is the hitting time on $A$. Harris-recurrence trivially implies $\phi$-irreducibility. A Harris-recurrent right process possesses an invariant measure [10]. In fact, when the invariant measure is finite, $X$ is called *positive Harris-recurrent*. A $\phi$-irreducible process is *aperiodic* if there exist an accessible $\nu_{\delta_m}$-petite set $C$ and $t_0$ such that $P^t(x, C) > 0$ for all $x \in C, t \geq t_0$. Meyn and Tweedie ([22], Proposition 6.1) show that a positive Harris-recurrent process is aperiodic if some skeleton chain $P^m$, $m > 0$, is irreducible, that is, if there exists a $\sigma$-finite measure $\phi$ on $\mathcal{B}(\mathcal{X})$ such that $\phi(A) > 0 \Longrightarrow \mathbb{E}_x[\sum_{n \geq 0} \mathbb{1}_A(X_{nm})] > 0$ for all $x \in \mathcal{X}$.

For Borel functions $f \geq 1$, $g$, define the norm $|g|_f := \sup_x |g(x)|/f(x)$ and the Banach space $\mathcal{L}_f := \{g, |g|_f < \infty\}$. For a signed measure $\mu$, the *total variation norm* is given by $\|\mu\|_{\mathrm{TV}} := \sup_A \mu(A) - \inf_A \mu(A)$; and the $f$-norm (for some Borel function $f \geq 1$), $\|\mu\|_f := \sup_{\{g, |g|_f = 1\}} |\mu(g)|$, so that the total variation norm is the $\mathbb{1}$-norm, where $\mathbb{1}$ denotes the constant function $\mathbb{1}(t) = 1$. The process is *ergodic* if

$$\forall x \in \mathcal{X}, \qquad \lim_{t \to \infty} \|P^t(x, \cdot) - \pi(\cdot)\|_{\mathrm{TV}} = 0,$$

and $f$-ergodic if $\pi(f) < \infty$ and

$$(1) \qquad \forall x \in \mathcal{X}, \qquad \lim_{t \to \infty} \|P^t(x, \cdot) - \pi(\cdot)\|_f = 0.$$

Finally, $X$ is *geometrically* (resp. *subgeometrically*) $f$-ergodic if the limit (1) holds at a rate $r(t) := \kappa^t$, for some $\kappa > 1$ [resp. $r := (r(t))_{t \geq 0}$ for some subgeometrically increasing rate]. A subgeometric rate is defined as follows (see, e.g., [34]). Let $\Lambda_0$ be the set of the measurable, bounded on bounded intervals and nondecreasing functions $r: \mathbb{R}_+ \to [1, \infty)$, such that $\log r(t)/t \downarrow 0$ as $t \to +\infty$. Let $\Lambda$ be the set of the rates $\bar{r} = (\bar{r}(t))_{t \geq 0}$ such that for some $r \in \Lambda_0$, $0 < \liminf_t \bar{r}(t)/r(t) \leq \limsup_t \bar{r}(t)/r(t) < \infty$. $\Lambda$ is by definition, the set of the subgeometric rates. For example, $\Lambda$ contains rates such as $\bar r(t) \sim \log^\beta(t+1)$, $\beta \geq 0$, $\bar r(t) \sim (1 \vee t^\alpha) \log^\beta(t+1)$, $\alpha > 0$, $\beta \in \mathbb{R}$, and subexponential rates $\bar r(t) \sim \exp(\alpha t^\beta)$, $\alpha > 0$, $0 < \beta < 1$.

Throughout this paper, we will often make comparison with (discrete time) Markov chains; the unfamiliar reader can refer to [21].

### 1.2. $(f,r)$-modulated moments and stability.

Define the hitting-time on a measurable set $C$, delayed by $\delta > 0$,

$$(2) \qquad \tau_C(\delta) := \inf\{t \geq \delta, X_t \in C\},$$

the moment $\tau_C(0)$ is denoted by $\tau_C$. It is proved in the literature that modulated moments of $\tau_C(\delta)$ for some closed petite set are related to Harris-recurrence, positive Harris-recurrence, $f$-ergodicity and geometric $f$-ergodicity. For a Borel function $f \geq 1$, an increasing nonnegative rate function $r = (r(t))_{t\geq 0}$, $\delta > 0$, define the $(f,r)$-modulated moment

$$G_C(x,f,r;\delta) := \mathbb{E}_x\left[\int_0^{\tau_C(\delta)} r(s)f(X_s)\,ds\right].$$

R1. $X$ is Harris-recurrent if and only if there exists a petite set $C$ such that, for all $x \in \mathcal{X}$, $\mathbb{P}_x(\tau_C < \infty) = 1$ ([20], Theorem 1.1).

R2. If $X$ is Harris-recurrent with invariant measure $\pi$, then for $f \geq 1$, $\pi(f) < \infty$ if and only if $\sup_{x \in C} G_C(x,f,\mathbb{1};\delta) < \infty$ for some closed petite set $C$ ([20], Theorem 1.2).

R3. A positive Harris-recurrent process is ergodic if and only if some skeleton chain $P^m$ is irreducible ([22], Theorem 6.1).

R4. A positive Harris-recurrent process is $f$-ergodic if (a) some skeleton chain $P^m$ is irreducible, (b) $\sup_{x \in C} G_C(x,h,\mathbb{1};\delta) < \infty$, where $h \geq \sup_{s \leq m} P^s f$ and $C$ is a closed petite set, and (c) for all $x$, $G_C(x,h,\mathbb{1};\delta) < \infty$ ([20], Proposition 4.1 and [22], Theorem 7.2).

R5. A positive Harris-recurrent process is geometrically $f$-ergodic if (a) some skeleton chain $P^m$ is irreducible, (b) there exists a closed petite set $C$ and $\eta > 0$ and $G_C(x,h,\exp(\eta t);\delta)$ is finite for all $x$, where $h \geq 1$ is a Borel function such that $c_1 f \leq \int_0^\infty \exp(-t)P^t h\,dt \leq c_2 f$ for some finite positive constants $c_i$, (c) $\sup_C G_C(\cdot,h,\exp(\eta t);\delta) < \infty$ ([5], Theorem 7.4).

In Section 2 we give a criterion of the form R1–R5 that implies subgeometric $f$-ergodicity.

To date, little is known about general characterizations for $f$-ergodicity at a subgeometric rate for Markov processes. However, we note some important special cases which have been studied in the literature. The work by Ganidis, Roynette and Simonot [9] is restricted (a) to convergence in total variation norm and (b) to diffusion processes on $\mathbb{R}^d$ with diffusion matrix equals to identity. Their proof is based on spectral properties of the transition semigroup seen as an operator, and differs from the probabilistic approach adopted in the present paper. We will see in Section 3.2 how to improve their conclusions. Veretennikov [36] and Malyshkin [19] deal with diffusion processes and can be read as a special application of the present paper. The most related work to the present one is the paper by Dai and Meyn [2] that considers $f$-ergodicity at a polynomial rate of a Markovian state process, in order to study the stability of multi-class queuing networks. These results are particularly related to our work and we will describe their results in Section 2.

1.3. *Drift condition and generator.* For a Borel function $0 \leq V < \infty$, denote by $\mathcal{A}V$ the Borel function—when exists—such that $t \mapsto \mathcal{A}V(X_t)$ is integrable $\mathbb{P}_x$-almost surely (a.s.), and there exists an increasing sequence of stopping-time $\{T_n\}_n$ such that for any stopping time $\tau$,

$$\mathbb{E}_x\left[V(X_{\tau \wedge T_n}) - V(X_0) - \int_0^{\tau \wedge T_n} \mathcal{A}V(X_s)\,ds\right] = 0 \qquad \text{for all } x \in \mathcal{X}, n \geq 0.$$

When $\mathcal{A}V$ exists, $V$ is said to be in the domain of $\mathcal{A}$. If there exists $h$ such that $t \mapsto h(X_t)$ is integrable $\mathbb{P}_x$-a.s. and $t \mapsto V(X_t) - V(X_0) - \int_0^t h(X_s)\,ds$ is a right-continuous $\mathbb{P}_x$-local martingale (with respect to the filtration $\mathcal{F}_t$), then $V$ is in the domain of $\mathcal{A}$ and $\mathcal{A}V = h$ [3]. If $V$ is in the domain of the weak infinitesimal generator $\tilde{\mathcal{A}}$, then $V$ is in the domain of $\mathcal{A}$ and $\mathcal{A}V = \tilde{\mathcal{A}}V$ [6]. If the functions $V$ and $\mathcal{A}V$ are right-continuous, these two sufficient conditions are equivalent and $\tilde{\mathcal{A}}V = \mathcal{A}V$.

When $\mathcal{A}V$ satisfies a drift condition $\mathcal{A}V \leq -f + b\mathbb{1}_C$ for some closed set $C$, and a nonnegative function $f$ such that $t \mapsto f(X_t)$ is integrable $\mathbb{P}_x$-a.s., we have $G_C(x, f, \mathbb{1}; \delta) \leq V(x) + \delta b$; this will be the basic tool to upper bound the $(f, r)$-modulated moments.

Conditions on $\mathcal{A}V$ are analogous to conditions on the variation $P^m V - V$ for a discrete time Markov chain with transition kernel $P^m$. It is well known that the condition $P^m V - V \leq -f$ outside a "test set" for the skeleton $P^m$ is related (a) to the $f$-ergodicity of the Markov chain $(X_{km})_k$ [21]; (b) to the geometric $V$-ergodicity if $f = \lambda V$ for some $0 < \lambda < 1$ [21]; (c) to the polynomial $V$-ergodicity if $f \propto V^{1-\alpha}$ for some $0 < \alpha \leq 1$ [8, 14]; (d) and more generally subgeometric $f$-ergodicity is $f \propto \phi(V)$ for some concave function $\phi$ [4]. Similar results hold for a continuous Markov process. Meyn and Tweedie [23] prove that the condition $\mathcal{A}V \leq -f$ outside a closed petite set is related (a) to the $f$-ergodicity of the Markov process $X$; and (b) to the geometric $V$-ergodicity if $f \propto V$ (see also [5, 26, 29]). In Section 2 we establish that the case $f \propto V^{1-\alpha}$ is related to polynomial $f$-ergodicity.

**2. Statements of the results.** In Theorem 1, we establish that modulated moment on some delayed hitting-time on a closed petite set $C$ provides a criterion for subgeometric $f$-ergodicity. We assume that there exist $\delta > 0$, a Borel function $f_* \geq 1$ and a rate function $r_* \in \Lambda$ such that

$$(3) \qquad \sup_C G_C(\cdot, f_*, \mathbb{1}; \delta) < \infty, \qquad \sup_C G_C(\cdot, \mathbb{1}, r_*; \delta) < \infty.$$

We will establish that $r_*$ is the maximal rate of convergence (that can be deduced from these assumptions) and it is associated to convergence in total variation norm, that is, in $\mathbb{1}$-norm, which is the minimal one. On the other hand, we will show that $f_*$ is the largest norm in which convergence occurs and the associated convergence rate is the minimal one $\mathbb{1}$.

Using an interpolation technique, we also derive a convergence rate $1 \leq r_f \leq r_*$ in $\Lambda$ associated to some $f$-norm, $1 \leq f \leq f_*$ (see [4] for a similar approach in the discrete time case). The simplest interpolation technique is given by Hölder's inequality which yields [from (3)] $\sup_C G_C(\cdot, f_*^p, r_*^{1-p}; \delta) < \infty$. By analogy to the discrete-time case, one would expect convergence in $f_*^p$-norm at the rate $r_*^{1-p}$, and we will prove the continuous time version of this result.

More generally, if there exists a pair of nondecreasing positive functions $(\Psi_1, \Psi_2)$ satisfying

$$(4) \qquad \Psi_1(x)\Psi_2(y) \leq x + y, \qquad x, y \geq 1,$$

then $\sup_C G_C(\cdot, \Psi_2(f_*), \Psi_1(r_*); \delta) < \infty$. We will establish that if $\Psi_1(r_*) \in \Lambda$, this condition yields convergence in $\Psi_2(f_*)$-norm at the rate $\Psi_1(r_*)$. Young functions are closely related to these pairs of functions $(\Psi_1, \Psi_2)$, say. Specifically, if $(H_1, H_2)$ is a pair of Young functions, then $(H_1^{-1}, H_2^{-1})$ satisfies (4) (see, e.g., [16], Chapter 1). Let $\mathcal{I}$ be the set of pairs of inverse Young functions, augmented with the pairs $(\mathrm{Id}, \mathbb{1})$ and $(\mathbb{1}, \mathrm{Id})$. As commented above, $\mathcal{I}$ contains the pairs $((x/p)^p, (y/(1-p))^{1-p})$, $0 < p < 1$, and, more generally, the pairs of functions increasing at infinity as $(x^p \ln^b x, y^{1-p} \ln^{-b} y)$ for some $0 < p < 1$ and $b \in \mathbb{R}$, $p = 0$ and $b \geq 0$, $p = 1$ and $b \leq 0$.

THEOREM 1. *Let $f_* \geq 1$ be a Borel function and $r_* \in \Lambda$. Assume that:*

(i) *$X$ is Harris-recurrent with invariant measure $\pi$, and some skeleton chain, say $P^m$, is $\psi$-irreducible.*

(ii) *There exist a closed petite set $C$ and some $\delta > 0$, such that (3) holds.*

(iii) *There exists a finite constant $c$ such that $\sup_{t \leq m} P^t f_* \leq c f_*$.*

*Then $\pi$ is an invariant probability measure, $\pi(f_*) < \infty$ and for any pair $\Psi := (\Psi_1, \Psi_2) \in \mathcal{I}$,*

$$\lim_{t \to +\infty} \{\Psi_1(r_*(t)) \vee 1\} \|P^t(x, \cdot) - \pi(\cdot)\|_{\Psi_2(f_*) \vee 1} = 0 \qquad \text{for all } x \in \mathcal{S}_\Psi,$$

where $\mathcal{S}_\Psi$, which is of $\pi$-measure one, is defined by

$$\mathcal{S}_\Psi := \left\{ x \in \mathcal{X}, \mathbb{E}_x \left[ \int_0^{\tau_C(\delta)} \Psi_1(r_*(s)) \Psi_2(f_*(X_s)) \, ds \right] < \infty \right\}.$$

The proof of Theorem 1 is postponed to Section 4.1. We first verify that $\Psi_1(r_*(t)) \vee 1 \in \Lambda$. Under (i) and (ii), $C$ is accessible and the following lemma holds ([22], Proposition 6.1).

LEMMA 2. *Suppose that $X$ is positive Harris-recurrent with invariant distribution $\pi$ and some skeleton chain $P^m$ is irreducible. For any accessible petite set $C$, there exist $t_0 > 0$ and an irreducibility measure $\nu$ for the process such that $\nu(C) > 0$ and $\inf_{x \in C} \inf_{t \geq t_0} P^t(x, \cdot) \geq \nu(\cdot)$.*

Based on this lemma and on condition (ii), the second step consists in proving that the skeleton $P^m$ is irreducible, aperiodic and possesses a petite set $A$ such that

$$\sup_{x \in A} G_A^{(m)}(x, \Psi_2(f_*), \Psi_1(r_*)) < \infty$$

$$(5)$$

$$\text{with } G_A^{(m)}(x, f, r) = \mathbb{E}_x \left[ \sum_{k=0}^{T_{m,A}} r(k) f(X_{km}) \right],$$

where $T_{m,A} \geq 1$ is the return-time to $A$ for the skeleton chain $P^m$

$$(6) \qquad T_{m,A} := \inf\{k \geq 1, X_{km} \in A\}.$$

By application of [35], Theorem 2.1, this proves that for $\pi$-a.a. $x$, $\lim_{k \to +\infty} \Psi_1(r_*(k)) \| P^{km}(x, \cdot) - \pi(\cdot) \|_{\Psi_2(f_*)} = 0$. Using Theorem 1(iii), the limit still holds replacing $r_*(k)$ (resp. $P^{km}$) by $r_*(t)$ (resp. $P^t$). We finally establish that the limit holds for all the points $x \in \mathcal{S}_\Psi$ and $\pi(\mathcal{S}_\Psi) = 1$.

REMARK 3. Theorem 1 remains valid by substituting condition (i) for the following condition: there exist a $\psi$-irreducible, aperiodic and positive recurrent transition kernel $P^m$.

Theorem 1 remains valid by substituting (ii) and (iii) for the following condition: there exist a closed petite set $C$ and some $\delta > 0$ such that $\sup_C G_C(\cdot, h, \mathbb{1}; \delta) < \infty$ and $\sup_C G_C(\cdot, \mathbb{1}, r_*; \delta) < \infty$ where $h \geq \sup_{t \leq m} P^t f_*$.

Condition (iii) implies that the semigroup $(P^t)_{t \geq 0}$ and the resolvent kernel $R = \int_0^\infty \exp(-t) P^t dt$ are bounded on $\mathcal{L}_{f_*}$.

REMARK 4. By (4), it is readily seen that $\{x; G_C(x, f_*, \mathbb{1}; \delta) + G_C(x, \mathbb{1}, r_*; \delta) < \infty\} \subset \mathcal{S}_\Psi$. It may be read from the proof that

$$(7) \qquad \lim_{t \to +\infty} \{\Psi_1(r_*(t)) \vee 1\} \| \mu P^t(\cdot) - \pi(\cdot) \|_{\Psi_2(f_*) \vee 1} = 0,$$

for all probability measure $\mu$ such that $G_C(x, \Psi_2(f_*) \vee 1, \Psi_1(r_*) \vee 1; \delta)$ is $\mu$-integrable. Applying again (4), (7) holds for all distribution $\mu$ such that $\{G_C(x, f_*, \mathbb{1}; \delta) + G_C(x, \mathbb{1}, r_*; \delta)\}$ is $\mu$-integrable.

REMARK 5. For any pair $(\Psi_1, \Psi_2) \in \mathcal{I}$, if $\Psi_1$ strongly increases at infinity [e.g., $\Psi_1(x) \propto x^p$ for some $p < 1$ close to one], then $\Psi_2$ slowly increases [$\Psi_2(x) \propto x^{1-p}$ for some $1 - p$ close to zero] ([16], Theorem 2.1, Chapter 1). Hence, the stronger the norm, the weaker the rate (and conversely). This compromise between the rate function and the norm of convergence is well known for the discrete parameter Markov chain ([35]; see also [4, 8, 14]). As expected, this property remains valid for the continuous-time Markov process.

Corollary 6 provides a condition based on $\mathcal{A}$, well adapted to prove polynomial ergodicity.

COROLLARY 6. Let $1 \leq V < \infty$ be a Borel function and $0 < \alpha \leq 1$. Assume that:

(i) Some skeleton chain $P^m$ is irreducible.

(ii) There exists a closed petite set $C$ such that $\sup_C V < \infty$ and for all $\alpha \leq \eta \leq 1$, $t \mapsto V^{\eta - \alpha}(X_t)$ is integrable $\mathbb{P}_x$-a.s. and

$$(8) \qquad \mathcal{A}V^\eta \leq -c_\eta V^{\eta - \alpha} + b \mathbb{1}_C, \qquad 0 \leq b < \infty, 0 < c_\eta < \infty.$$

Then there exists a unique invariant distribution $\pi$, $\pi(V^{1-\alpha}) < \infty$ and for all $0 < p < 1$ and $b \in \mathbb{R}$ or $p = 1$ and $b \geq 0$ or $p = 0$ and $b \leq 0$,

$$\lim_{t \to +\infty} (1 + t)^{(1-p)(1-\alpha)/\alpha} (\log t)^b \|P^t(x, \cdot) - \pi(\cdot)\|_{V^{(1-\alpha)p}(\ln V)^{-b} \vee 1} = 0,$$

$$x \in \mathcal{X}.$$

The proof is given in Section 4.2. From (ii), we obtain $G_C(x, V^{1-\alpha}, \mathbb{1}; \delta) + G_C(x, \mathbb{1}, (1+t)^{1/\alpha-1}; \delta) \leq cV(x)$; and then we apply Theorem 1.

By choosing $b = 0$ and $p = (1 - \kappa\alpha)/(1 - \alpha)$ for some $1 \leq \kappa \leq 1/\alpha$, Corollary 6 yields

$$(9) \qquad \forall x \in \mathcal{X}, \qquad \lim_{t \to \infty} (t+1)^{\kappa - 1} \|P^t(x, \cdot) - \pi(\cdot)\|_{V^{1-\kappa\alpha}} = 0.$$

If (9) holds for some $V$ function, we shall say that the Markov chain is *polynomially ergodic* with rate $(1 + t)^{(1-\alpha)/\alpha}$.

REMARK 7. Corollary 6 can be compared to the paper by Jarner and Roberts [14] for the discrete parameter case. They start with proving that if there exist a Borel function $1 \leq V < \infty$, $0 < \alpha \leq 1$, a set $C$ such that, for all $\alpha \leq \eta \leq 1$,

$$(10) \qquad P^m V^\eta - V^\eta \leq -c_\eta V^{\eta-\alpha} + b\mathbb{1}_C, \qquad 0 \leq b < \infty, 0 < c_\eta < \infty,$$

there exists $c < \infty$ such that $G_C^{(m)}(\cdot, V^{1-\alpha}, \mathbb{1}) + G_C^{(m)}(\cdot, \mathbb{1}, (1+t)^{1/\alpha-1}) \leq cV$, where $G_C^{(m)}$ is given by (5). The drift condition (10) is analogous to (8) and the controls of the moments $G_C^{(m)}$ and $G_C$ are similar. If, in addition, $P^m$ is irreducible, aperiodic and $C$ is petite for the skeleton, $P^m$ is positive with invariant distribution $\pi$ such that $\pi(V^{1-\alpha}) < \infty$ and for all $1 \leq \kappa \leq 1/\alpha$, the skeleton is $V^{1-\kappa\alpha}$-ergodic with rate $(n+1)^{\kappa-1}$. These rates coincide with those in (9).

REMARK 8. From the proof of Corollary 6, it may be read than only a finite number of nested drift conditions is required; nevertheless, in practice, it is not more restricting to verify a continuum of drift conditions than to verify a finite number of drift conditions. More precisely, assumption (ii) can be substituted for the following conditions: (iii) there exist a closed petite set and functions $1 \leq V_{q-1} \leq cf_q$, such that for all integers $1 \leq q \leq p$, $\mathcal{A}V_q \leq -f_q + b\mathbb{1}_C$, $t \mapsto f_q(X_t)$ is integrable $\mathbb{P}_x$-a.s., and $\sup_C V_p < \infty$; (iv) there exists $\beta > 0$ such that $\mathbb{E}_x[\tau_C^\beta] \leq f_1(x)$. If such, following the same lines as in the proof of Proposition 26, it may be proved that $G_C(\cdot, \mathbb{1}, (t+1)^{p-1+\beta}; \delta) + G_C(\cdot, f_p, \mathbb{1}; \delta) \leq cV_p$ for some $\delta > 0$. Together with condition (i), this yields $f_*^{1-\eta}$-ergodicity at a rate $(t+1)^{(p-1+\beta)\eta}$ for all $0 \leq \eta \leq 1$, where $f_*$ is any function satisfying $\sup_{t \leq m} P^t f_* \leq f_p$.

Dai and Meyn [2] are, to our best knowledge, the first to exhibit this kind of nested drift condition and, hence, the first to address ergodicity at a polynomial rate; they proved this yields $f_1$-ergodicity at a rate $(t+1)^{p-1}$ (Theorem 6.3, [2]). We are able to obtain the same result: to that goal, we observe that conditions (iii) and (iv) are verified with functions $f_k \leq f_p^{k/p}$, $\beta = 1$ (as a consequence of Proposition 5.3 and equation (6.1) in [2]) and $f_* \propto f_p$.

We proved that nested drift conditions on the generator $\mathcal{A}$ provide a control of moments $G_C$ with a polynomially increasing rate function. The converse seems to be an open question. We, nevertheless, make mention of Propositions 5.4 and 6.1 in [2], that provide a (partial) converse condition: from the condition $\sup_C G_C(\cdot, f, \mathbb{1}; \delta) < \infty$, they deduce a drift condition on $\mathcal{A}$ (we point out that this single condition implies a continuum of conditions by using the same convexity argument as in [14], Lemma 3.5). Unfortunately, this drift condition, in turn, implies only a control of the moment $G_C(x, Rf, \mathbb{1}; \delta)$, where $Rf(x)$ is a function, which is, in general, difficult to compare with $f$.

**3. Examples.** In this section $\langle \cdot, \cdot \rangle$ and $|\cdot|$ denote, respectively, the scalar product and the Euclidean norm in $\mathbb{R}^n$. If $u$ is a twice continuously differentiable real-valued function on $\mathbb{R}^n$, $\nabla u$ (resp. $\nabla^2 u$) denotes its gradient (resp. its Hessian matrix); and $\partial u / \partial x_i$ its partial derivative with respect to the $i$th variable. For a matrix $u$, $\mathrm{Tr}(u)$ stands for the matrix trace and $u'$ the matrix transpose. For $r \in \Lambda$, define the sequence $r^0$ by $r^0(t) = \int_0^t r(s) \, ds$. Finally, we largely make use of the inequality $r(s+t) \le r(s)r(t)$, $s, t \ge 0$, which holds for any rate $r \in \Lambda$ ([34], Lemma 1).

3.1. *Toy example*: *jump process*. Consider the jump process on $\mathbb{Z}_+$ such that given that $X_t = i$, the waiting-time to the next jump has an exponential distribution with expectation $\lambda_i^{-1}$ and is independent of the past history. We assume that for all $i \ge 0$, $\lambda_i > 0$, and $\sup_{i \ge 0} \lambda_i < \infty$. The probability that the jump leads to state $j$ is given by the matrix entry $Q(i, j)$. We consider the case when $Q(0, i) = p_i$ and $Q(i, 0) = 1$ for all $i \ge 1$, for some positive sequence $(p_i)_{i \ge 1}$ such that $\sum_{i \ge 1} p_i = 1$. We assume, in addition, that

$$(11) \qquad \liminf_i \lambda_i = 0 \quad \text{and} \quad \sum_{i \ge 0} p_i \lambda_i^{-1} < \infty.$$

Since $\sup_{i \ge 0} \lambda_i < \infty$, there exists a $\mathbb{Z}_+$-valued right-continuous strong Markov process satisfying the heuristic description above and such that, for all $(i, j) \in \mathbb{Z}_+^2$, the limit exists

$$(12) \qquad \lim_{t \to 0} \frac{P^t(i, j) - \delta_i(j)}{t} =: \mathsf{A}(i, j) < \infty,$$

where $\delta_i$ is the Dirac-mass at point $i$, and for all $i \ge 1$,

$$(13) \qquad \mathsf{A}(0, 0) = -\lambda_0, \qquad \mathsf{A}(0, i) = \lambda_0 p_i, \qquad \mathsf{A}(i, 0) = -\mathsf{A}(i, i) = \lambda_i,$$

and $\mathsf{A}(i, j) = 0$ otherwise (see, e.g., [7], page 330).

LEMMA 9. *The process is Harris-recurrent, reversible with invariant distribution $\pi$ given by $\pi(0) = \{1 + \sum_{j \ge 1} p_j \lambda_j^{-1}\}^{-1}$ and $\pi(i) = p_i \lambda_i^{-1} \pi(0)$, $i \ge 1$. Any skeleton chain $P^m$ is irreducible.*

PROOF. We have $\mathbb{E}_i[\tau_0] = (1 - \mathbb{1}_0(i)) \lambda_i^{-1}$ and for all $i \ge 0$, $j \ne i$, $\delta > 0$, $P^\delta(i, j) \le p_j$. Then, $\mathbb{E}_i[\tau_0(\delta)] = \delta + \sum_{j \ge 1} P^\delta(i, j) \mathbb{E}_j[\tau_0] \le \delta + 2 \sum_{j \ge 1} p_j \lambda_j^{-1}$. Hence, for all $i \in \mathbb{Z}_+$, $\mathbb{P}_i(\tau_0(\delta) < \infty) = 1$ and as $\{0\}$ is a closed petite set, the process is Harris-recurrent. $\pi$ is the unique invariant probability measure (as unique solution to $\pi \mathsf{A} = 0$), and since $X$ obeys the detailed balance, that is, $\pi(i) \mathsf{A}(i, j) = \pi(j) \mathsf{A}(j, i)$ for all $i, j$, the process is reversible. Finally, for all $m > 0$, and $i, j > 0$,

$$P^m(i, j) \ge p_j \lambda_0 \lambda_i \lambda_j \int_0^m ds \exp(-\lambda_i s)$$

$$\times \int_0^{m-s} dt \exp(-\lambda_0 t) \int_0^{m-(t+s)} du \exp(-\lambda_j u) > 0,$$

where the inequality says that $P^m(i,j)$ is greater than the probability of a single visit to 0 before a jump to $j$. Similarly, it is easy to prove that $P^m(0,j) > 0$ and $P^m(j,0) > 0$ for any $j \in \mathbb{Z}_+$. This proves the irreducibility of any skeleton. $\square$

We deduce from R3 that the process is ergodic. Nevertheless, this convergence fails to occur at a geometric rate as shown in Lemma 10, the proof of which relies on the notion of conductance.

LEMMA 10. *$X$ fails to be geometrically ergodic.*

PROOF. As $X$ is reversible, any Markov kernel $P^m$ is reversible. It is proved in [17] that for a reversible Markov kernel $P^m$, the conductance $\kappa_m$ given by $\kappa_m := \inf_A c_m(A)$, where $c_m(A) := \{\pi(A)\pi(A^c)\}^{-1} \int_A P^m(x, A^c)\pi(dx)$, is positive if and only if $P^m$ is geometrically ergodic. We verify that for any skeleton $P^m$, the conductance is zero which will involve that the skeleton fails to be geometrically ergodic. Consider the set of states $i$ such that $\pi(i) \le 1/2$. Then $c_m(i) \le 2(1 - P^m(i,i)) \le 2(1 - \exp(-\lambda_i m))$ upon noting that $P^m(i,i)$ is lower bounded by the probability that the waiting-time in state $i$ is greater than $m$. Since $\liminf_{i \to +\infty} \lambda_i = 0$, for all $\varepsilon > 0$, there exists a state $i$ such that $c_m(i) \le \varepsilon$, which involves $\kappa_m = 0$. $\square$

We now identify functions $V$ that are in the domain of $\mathcal{A}$.

LEMMA 11. *Let $0 \le V < \infty$ be a Borel function such that $\sum_{i \ge 1} p_i V(i) < \infty$. Then $V$ is in the domain of $\mathcal{A}$ and $\mathcal{A}V = \mathsf{A}V$.*

PROOF. For a function $f \ge 0$ such that $\sum_j p_j f(j) < \infty$, the monotone convergence theorem yields

$$\sum_{j \ge 1} \lim_{t \downarrow 0} \frac{P^t(i,j) - \delta_i(j)}{t} f(j) = \lim_{J \uparrow +\infty} \sum_{j=0}^{J} \lim_{t \downarrow 0} \frac{P^t(i,j) - \delta_i(j)}{t} f(j) = \mathsf{A}f(i);$$

in addition, $\sup_{\mathbb{Z}_+} f^{-1}|\mathsf{A}f| < \infty$. This proves that $V$ is in the domain of the weak infinitesimal generator $\mathsf{A}$, and thus in the domain of $\mathcal{A}$. $\square$

The expression of the generator suggests that function $V$ on the form $\lambda_i^{-\rho}$ is a candidate to solve the drift inequality (8). This yields $f$-ergodicity at a log-polynomial rate.

PROPOSITION 12. *Assume that there exists $\beta \ge 1$ such that $\sum_{i \ge 1} p_i \lambda_i^{-\beta} < \infty$. Then for all $i \in \mathbb{Z}_+$, $0 < \kappa < \beta - 1$ and $b \in \mathbb{R}$ or $\kappa = 0$ and $b \le 0$, or $\kappa = \beta - 1$ and $b \ge 0$,*

$$\lim_{t \to \infty} (1+t)^{\beta-1-\kappa} [\ln(1+t)]^b \|P^t(i, \cdot) - \pi(\cdot)\|_{1 + \lambda_x^{-\kappa}[\ln(1 + \lambda_x^{-1})]^{-b}} = 0.$$

PROOF. We apply Corollary 6: we choose $V \geq 1$ such that for all $i \geq 1$, $V(i) = c^{-1}V(0)\lambda_i^{-\beta}$ for some $c > V(0)$. Then (8) is verified with $\alpha = \beta^{-1}$ and the closed petite set $C = \{0\}$. □

When $\beta = 1$ [i.e., with nothing more than the condition (11)], this establishes the convergence in total variation norm at the rate $\mathbb{1}$, which corroborates the ergodicity of the process proved above. Nevertheless, if for some $\beta > 0$, the sum $\sum_{i \geq 1} p_i (1 \vee \lambda_i^{-1})[\log(1 \vee \lambda_i^{-1})]^\beta$ exists, Corollary 6 does not yield a stronger convergence result than the ergodic one. We prove, by application of Theorem 1, that covers more general rates than the polynomial ones, that convergence in total variation norm occurs at the rate $r_*(t) \sim [\log(t)]^\beta$, and convergence in norm $f_*(x) = [\log(1 \vee \lambda_x^{-1}) + 1]^\beta$ occurs at rate $\mathbb{1}$. We also derive sufficient conditions for subexponential ergodicity.

LEMMA 13. Let $f_* : \mathbb{Z}_+ \to [1, \infty)$ and $r_* \in \Lambda$ such that

$$\sum_{j \geq 1} p_j (1 \vee \lambda_j^{-1}) f_*(j) < \infty \quad and$$

(14)

$$\sum_{j \geq 1} p_j \lambda_j^{-1} \int_0^{+\infty} r_*(s) \lambda_j \exp(-\lambda_j s)\, ds < \infty.$$

Then there exists a finite constant $c$ such that for all $m > 0$, $\sup_{t \leq m} P^t f_* \leq c f_*$. For all $\delta > 0$, there exists a finite constant $c$ such that

$$G_0(x, f_*, \mathbb{1}; \delta) \leq c(1 \vee \lambda_x^{-1}) f_*(x),$$

$$G_0(x, \mathbb{1}, r_*; \delta) \leq c \int_0^{+\infty} r_*(s) \exp(-\lambda_x s)\, ds.$$

PROOF. Since $P^t(x, j) \leq p_j$ for all $x \neq j$, it is trivial to prove that $\sup_{t > 0} \sup_{i \in \mathbb{Z}_+} f_*^{-1} P^t f_* < \infty$. For $f \geq 1$ and $r \in \Lambda$,

$$G_0(x, f, r; \delta) \leq \int_0^\delta r(s) P^s f(x)\, ds + r(\delta) \sum_{j \geq 1} P^\delta(x, j) f(j) \mathbb{E}_j[r^0(\tau_0)].$$

$\mathbb{E}_j[r^0(\tau_0)] = \lambda_j \int r^0(t) \exp(-\lambda_j s)\, ds = \int r(t) \exp(-\lambda_j s)\, ds$. □

PROPOSITION 14. (i) Assume that $\sum_{i \geq 1} p_i (1 \vee \lambda_i^{-1})[\log(1 \vee \lambda_i^{-1})]^\beta < \infty$, for some $\beta \geq 0$. For all $0 \leq \kappa \leq \beta$, $i \in \mathbb{Z}_+$, $\lim_{t \to +\infty} [\log(t+1)]^{\beta - \kappa} \|P^t(i, \cdot) - \pi(\cdot)\|_{[1 + \log(1 \vee \lambda_i^{-1})]^\kappa} = 0$.

(ii) Assume that $\sum_{i \geq 1} p_i (1 \vee \lambda_i^{-1}) \lambda_i^{-1/2} \exp(z^2 \lambda_i^{-1}) < \infty$, for some $z > 0$. For all $0 \leq p \leq 1$, $i \in \mathbb{Z}_+$, $\lim_{t \to +\infty} \exp(2z(1-p)t^{1/2}) \|P^t(i, \cdot) - \pi(\cdot)\|_{[1 + \lambda_i^{-1/2} \exp(z^2 \lambda_i^{-1})]^p} = 0$.

Proof.    In both cases, apply Theorem 1; for case (i), set $r_*(t) = \{\log(\exp(\beta - 1) + t\}^\beta$ and $f_*(i) = 1 + \log(1 \vee \lambda_i^{-1})^\beta$; and for case (ii), set $r_*(t) = \exp(2zt^{1/2})$, $f_*(i) = 1 + \lambda_i^{-1/2} \exp(z^2 \lambda_i^{-1})$ and observe that $\int \exp(2zs^{1/2}) \lambda \exp(-\lambda s) \, ds \leq 1 + 2\sqrt{\pi} z \lambda^{-1/2} \exp(z^2 \lambda^{-1})$.    □

3.2. *Langevin tempered diffusions on $\mathbb{R}^n$.*    Let us consider a stochastic integral equation

$$\text{(15)} \qquad X_t = X_0 + \int_0^t b(X_s) \, ds + \int_0^t \sigma(X_s) \, dW_s,$$

where $W_t$ is an $n$-dimensional Brownian motion, the drift coefficient $b = (b_1, \ldots, b_n)'$ is on the form, $1 \leq i \leq n$,

$$b_i(x) = \frac{1}{2} \sum_{j=1}^n a_{i,j}(x) \frac{\partial}{\partial x_j} \log \pi(x) + \frac{1}{2} \sum_{j=1}^n \frac{\partial}{\partial x_j} a_{i,j}(x),$$

where $a = \sigma\sigma'$ is the $n \times n$ symmetric positive definite matrix. Such a diffusion is the so-called Langevin diffusion and is defined in such a way that $\pi$ is, up to a multiplicative constant, the density of the unique invariant probability distribution (with respect to the Lebesgue measure on $\mathbb{R}^n$). This property motivates recent interests in Langevin diffusion for their use as MCMC methods, where the scope of these techniques is to draw samples from a Markov chain with given stationary density $\pi$. The efficiency of these algorithms is linked to the rate at which $f$-moments $\mathbb{E}_x[f(X_t)]$ converge to the constant $\pi(f)$. This motivates the study of the $f$-ergodicity. In practice, discretizations of the continuous-time process are used to solve the MCMC simulation problem and recent works proved that it is possible to find methods of discretizing which inherit the convergence rates of the continuous-time diffusion (see [27, 28, 30, 32, 33] for methods of discretizing and their use in MCMC techniques). Roberts and Tweedie proved that, on the real line, when the target density $\pi$ is heavy tailed, the Langevin diffusion with $a := 1$ cannot be geometrically ergodic. We complement this assertion when $\pi$ is polynomial in the tails, and prove that the Langevin diffusion in the one-dimensional case, as well as in the multidimensional one, is $f$-ergodic at a polynomial rate. For such polynomial target density on the real line, it was observed in [13] that the polynomial rate of convergence of the Metropolis–Hastings algorithm could be improved by choosing a heavy-tailed proposal distribution. This idea, when adapted to the diffusion on the real line, suggests the choice of a speed measure, that is, of the coefficient $\sigma$ such that $\sigma$ is small when the process is close to the modes of $\pi$ and big when far from the modes [32]. In the multidimensional case, this suggests $a(x)$ on the form $\pi^{-2d}(x)\mathbb{I}_n$, where $\mathbb{I}_n$ is the identity matrix on $\mathbb{R}^n$, $d > 0$. In that case $(d > 0)$, we call these processes Langevin tempered diffusion (see [27] for the justification of these heated diffusions). It was observed in the literature that by choosing $d$ large enough, a diffusion on the real line with target density polynomial in the tails is geometrically ergodic. We investigate the behavior of this Langevin tempered diffusion in the multidimensional case, contrarily to most of the mentioned contributions that cover the one-dimensional case. In Theorem 16, it is proved that, up to some critical temperature $d_*$, the diffusion is polynomially ergodic and that the larger $d$, the better the rate. When $d \geq d_*$, the diffusion is geometrically ergodic. We henceforth consider a diffusion matrix $a(x) = \sigma^2(x)\mathbb{I}_n$, where $\sigma(x) := \pi^{-d}(x)$ for some $d \geq 0$. Assume the following:

A1. $\pi$ is, up to a multiplicative constant, a positive and twice continuously differentiable density on $\mathbb{R}^n$ (with respect to the Lebesgue measure).

Define the drift vector

$$(16) \quad b(x) := \frac{1}{2}\sigma^2(x)(\nabla \log\{\pi(x)\sigma^2(x)\}) = \frac{1-2d}{2}\pi^{-2d}(x)\nabla\log\pi(x).$$

Under A1, the coefficients $b$ and $\sigma$ are locally Lipschitz-continuous, which implies that for any compact set $\mathcal{K}$, $\sup_{x \in \mathcal{K}}\{|b(x)| + |\sigma(x)|\}(1 + |x|)^{-1} < \infty$. These local conditions allow the construction of a continuous process satisfying the stochastic integral equation (15) up to the explosion time $\zeta := \lim_{n\to\infty}\zeta_n$, where $\zeta_n := \inf\{t \geq 0, |X_t| \geq n\}$. We thus formulate the following assumption:

A2. The process is regular, that is, $\zeta = +\infty$ a.s.

Under A1, a sufficient condition for regularity is the existence of a twice continuously differentiable nonnegative function $V$ and a constant $c \geq 0$ such that $\mathsf{L}V \leq cV$ on $\mathbb{R}^n$ and $\lim_{n\to\infty}\inf_{|x|\geq n}V(x) = +\infty$ ([12], Theorem 3.4.1), where $\mathsf{L}$ is the elliptic operator

$$\mathsf{L}V(x) = \langle b(x), \nabla V(x)\rangle + \frac{\mathrm{Tr}(\nabla^2 V(x)a(x))}{2}$$
$$= \frac{\pi^{-2d}(x)}{2}\left((1-2d)\langle\nabla\log\pi(x), \nabla V(x)\rangle + \sum_{i=1}^{n}\frac{\partial^2 V(x)}{\partial x_i^2}\right).$$

In the one-dimensional case, Has'minskii ([12], Remark 2, page 105) establishes that the process is regular if $d$ is chosen such that

$$(17) \qquad\qquad \int_{\mathbb{R}}\pi^{2d-1}(x)\,dx = +\infty,$$

since the function $V(x) := \mathrm{sign}(x)\int_0^x Q(y)\,dy$, where $\ln Q(x) = -2\int_0^x b(t) \times \sigma^{-2}(t)\,dt = (2d-1)(\ln\pi(x) - \ln\pi(0))$ is finite and satisfies $\mathsf{L}V = 0$ on $\mathbb{R}$. To cover the multidimensional case, we adapt this condition and claim that the process is regular if $d$ is chosen such that

$$(18) \quad \int_r^\infty t^{1-n} \exp\left(-(1-2d)\int_r^t s^{-1} \sup_{\{x,|x|=s\}} \langle \nabla \log \pi(x), x \rangle \, ds\right) dt = +\infty.$$

Indeed, the function $V(x) := U(|x|)$ where, for all $u \geq 0$,

$$U(u) := \int_r^u \exp\left(-\int_r^t \sup_{\{x,|x|=s\}} \left\{\left\langle \frac{2b(x)}{\sigma^2(x)}, \frac{x}{|x|} \right\rangle + \frac{n-1}{|x|} \right\} ds\right) dt,$$

is finite and satisfies $\mathsf{L}V = 0$ on $\mathbb{R}^n$. In the one-dimensional case, condition (17) is necessary for the existence of an invariant probability measure ([12], Remark 2, page 105); thus, for the objective of this paper, $d$ has to be chosen in the set $\mathcal{D}_1$ of the positive real numbers such that (17) hold. Observe that $\mathcal{D}_1$ is nonempty and contains $\{0, 1/2\}$. In the multidimensional case, a necessary condition for (positive) recurrence is that $d$ checks a condition on the form (18) where the supremum is replaced by the infimum ([11], Theorem II, page 194). This involves the definition of an interval $\mathcal{D}_n$ limiting the range of the possible temperature $d$.

Under A1 and A2, there exists a solution $(\Omega, \mathcal{F}, (\mathcal{F}_t), (W_t), (X_t), \mathbb{P})$, where $(\Omega, \mathcal{F}, (\mathcal{F}_t), (W_t), \mathbb{P})$ is $n$-dimensional Brownian motion, $(X_t)_t$ is an $\mathcal{F}_t$-adapted homogeneous and continuous Markov process with Feller transition probability, satisfying (15) $\mathbb{P}$-a.s. and such that both the integral exist, that is, for all $t > 0$,

$$(19) \quad \mathbb{P}\left(\int_0^t b(X_s)\,ds + \int_0^t \sigma^2(X_s)\,ds < \infty\right) = 1.$$

A transition semigroup $(P^t)_{t \geq 0}$ has the Feller property if for any continuous bounded real-valued function $f$, $x \mapsto P^t f(x)$ is continuous. $(X_t)$ is thus a strongly Markovian process as a (right)-continuous process with Feller transition probability [6].

Let $0 \leq V < \infty$ be a twice continuously differentiable function such that there exist a nonnegative Borel function $\phi$, bounded on compact sets, a constant $b < \infty$ and a compact set $C$ such that $\mathsf{L}V \leq -\phi \mathbb{1}_{C^c} + b\mathbb{1}_C$. From (19) and the continuity of $t \mapsto \nabla V(X_t)$, the process $t \mapsto \int_0^t \sigma(X_s)\{\nabla V(X_s)\}' \, dW_s$ is a local martingale. Application of Itô's rule yields $\mathsf{L}V = \mathcal{A}V$.

A3. For all $1 \leq i, j \leq n$, $\partial^2 \sigma^2(x)/\partial x_i \, \partial x_j$ and $\partial^2 \log \pi(x)/\partial x_i \, \partial x_j$ are locally uniformly Hölder continuous.

PROPOSITION 15. *Under* A1–A3, *the process is reversible and* $\pi$ *is, up to a multiplicative constant, the density of an invariant probability measure. Any skeleton chain is irreducible, and compact sets are closed petite sets.*

PROOF. There exists a continuous function $p : (t, x, y) \mapsto p(t, x, y)$ such that $P^t(x, dy) = p(t, x, y) \, dy$ ([15], Theorem 1.1). Since the process is regular (or conservative, in the terminology of Kent [15]) and $\pi$ is Lebesgue integrable ([15], Theorems 4.1. and 6.2) imply that the process is time-reversible and

$$(20) \qquad \lim_{t \to \infty} \int_A p(t, x, y) \, dy = \left( \int \pi(x) \, dx \right)^{-1} \int_A \pi(x) \, dx.$$

Hence, $\pi(dx)$ is invariant. Irreducibility of skeletons results from (20), and petiteness of compact sets from the continuity of $p(t, \cdot, \cdot)$. □

Finally, we restrict our attention to densities $\pi$ that are polynomially decreasing in the tails.

A4. $\pi$ satisfies A1 and A3 and there exists some $0 < \beta < 1/n$,

$$0 < \liminf_{|x| \to +\infty} \frac{|\nabla \log \pi(x)|}{\pi^\beta(x)} \leq \limsup_{|x| \to +\infty} \frac{|\nabla \log \pi(x)|}{\pi^\beta(x)} < \infty,$$

$$2\beta - 1 < \liminf_{|x| \to +\infty} \frac{\mathrm{Tr}(\nabla^2 \log \pi(x))}{|\nabla \log \pi(x)|^2} \leq \limsup_{|x| \to +\infty} \frac{\mathrm{Tr}(\nabla^2 \log \pi(x))}{|\nabla \log \pi(x)|^2} < \infty.$$

Set $\gamma = \liminf_{|x| \to +\infty} \mathrm{Tr}(\nabla^2 \log \pi(x)) |\nabla \log \pi(x)|^{-2}$.

This class is nonempty and contains the densities that are polynomially decreasing in the tails $\pi(x) = c|x|^{-1/\beta}$ for large $|x|$, where $0 < \beta < 1/n$; in that case, $\gamma = \beta(2 - n) > 2\beta - 1$. For this family, the regularity criterion (17) or (18) says that the temperature $d$ has to be chosen in $\mathcal{D}_n = [0, (1 + \beta(2-n))/2]$. For any density in the class A4, $0 < \liminf_{|x|} |x|^\beta \pi(x) \leq \limsup_{|x|} |x|^\beta \pi(x) < \infty$. Hence, $\mathcal{D}_1 = [0, (1+\beta)/2]$ and for $n \geq 3$, $1/2 \notin \mathcal{D}_n$. If $\sup_{s \geq r} \sup_{\{x, |x|=s\}} \langle \nabla \log \pi(x), x \rangle =: -\varrho^{-1} < 0$ exists, then $[0, 1/2 + \varrho(1 - n/2)] \subseteq \mathcal{D}_n$.

It is readily seen that setting $V = 1 + \mathrm{sign}(\rho)\pi^{-\rho}$ outside a compact set, and $V = 1$ otherwise,

$$(21) \qquad \mathsf{L}V = -\frac{|\rho|}{2} V \frac{\pi^{-\rho}}{1 + \pi^{-\rho}} \pi^{2(\beta-d)}$$
$$\times \left( \frac{|\nabla \log \pi|}{\pi^\beta} \right)^2 \left( 1 - \rho - 2d + \frac{\mathrm{Tr}(\nabla^2 \log \pi)}{|\nabla \log \pi|^2} \right),$$

for large $|x|$. As established in [32], Theorem 3.1, the diffusion cannot be geometrically ergodic when $0 \leq d < \beta$: by choosing $f := \pi^{d-\beta}$ and applying It's formula, $df(X_t) \approx c_1 \pi^{\beta-d}(X_t) \, dt + c_2 \, dW_t$ for some constants $c_i$; and the drift coefficient tends to zero for a large value of the process. The process $(f(X_t))_t$ fails to be geometrically ergodic, and, henceforth, $(X_t)_t$ itself.

From (21), for large $|x|$,

$$\mathsf{L}V \leq -cV^{1-\alpha} \qquad \text{where } \alpha := 2\rho^{-1}(\beta - d) \text{ and } c > 0$$

$$\iff \quad 1 + \gamma - \rho - 2d > 0.$$

In any cases, one has to choose $\rho$ such that $c > 0$. If $\alpha \leq 0$ and $\rho > 0$, then the process is geometrically $V$-ergodic ([23], Theorem 6.1). If $0 < \alpha \leq 1$ and $\rho > 0$, the diffusion is polynomially ergodic as discussed in Section 2. If $\alpha \leq 1$ and $\rho$ can be set negative, the process is uniformly ergodic, that is, there exist $\kappa > 1$ and a constant $c$ such that for all $x$, $\lim_{t \to \infty} \kappa^t \|P^t(x, \cdot) - \pi(\cdot)\|_{\mathrm{TV}} \leq c$ and the convergence does not depend on the starting point. This yields Theorem 16: the first assertion results from [28] and Corollary 6 of the present paper. The second and third assertions result from [23], Theorem 6.1.

THEOREM 16.  *Consider the Langevin tempered diffusion on $\mathbb{R}^n$, where the target density $\pi$ is from the class* A4 *and $\sigma := \pi^{-d}$ for some $d$ satisfying* (17) *if $n = 1$ or* (18) *if $n \geq 2$.*

(i) *If $0 \leq d < \beta$, the process fails to be geometrically ergodic. For all $0 \leq \kappa < 1 + \gamma - 2\beta$,*

$$(22) \qquad \lim_{t \to +\infty} (t+1)^\tau \|P^t(x, \cdot) - \pi(\cdot)\|_{1 + \pi^{-\kappa}} = 0, \qquad \tau < \frac{1 + \gamma - 2\beta - \kappa}{2(\beta - d)}.$$

(ii) *If $\beta \leq d < (1 + \gamma)/2$, then for all $0 < \kappa < 1 + \gamma - 2d$, the diffusion is geometrically $V$-ergodic with $V := 1 + \pi^{-\kappa}$.*

(iii) *If $\beta < d$, the diffusion is uniformly ergodic.*

Theorem 16 extends earlier results to the multidimensional case and provides polynomial rates of convergence of the "cold" Langevin tempered diffusions for a wide family of norms. In the one-dimensional case, when $d = 0$ ([9], Result R3, page 245) only claim that the convergence in total variation norm is polynomial, with no explicit value of the rate of convergence. We establish that, for a given $\pi^{-\kappa}$-norm, the minimal rate of convergence is achieved with $d = 0$ and, in that case, coincides with the rate of convergence of the symmetric random walk Hastings–Metropolis algorithm with light proposal distribution [13]. By choosing a diffusion matrix which is heavy where the target distribution is light, and conversely, improves the rate of convergence as evidenced by (22). The critic temperature is $d = \beta$. For $d \geq \beta$, the diffusion is no more polynomially ergodic and geometric rates can be reached. This critic temperature coincides with the critic one given in [32], Theorem 3.1, for the real-valued diffusion.

REMARK (*General diffusions on $\mathbb{R}^n$*).    The techniques above can be adapted for the analyzing of diffusions satisfying (15). Under conditions implying (a) the existence of a solution, (b) the condition (i) of Corollary 6 and (c) the petiteness property of the compact sets (see, e.g., [12, 19, 36]), we are able to prove that when there exist $M, \beta, \gamma > 0$ and $l < 2$ such that

$$\sup_{\{x, |x| \geq M\}} |x|^{-(2+l)} \langle x, a(x)x \rangle =: \beta, \qquad \sup_{\{x, |x| \geq M\}} |x|^{-l} \operatorname{Tr}(a(x)) =: \gamma,$$

$$\sup_{\{x, |x| \geq M\}} |x|^{-l} \langle b(x), x \rangle =: -r \qquad \text{for some } r > (\gamma - \beta l)/2,$$

then the diffusion is polynomially ergodic and for all $x$, for all $0 \leq \kappa < l + \beta^{-1}(2r - \gamma)$,

$$\lim_{t \to \infty} (1 + t)^\tau \| P^t(x, \cdot) - \pi(\cdot) \|_{1 + |x|^\kappa} = 0, \qquad \tau < \frac{2(r + \beta) - \gamma}{\beta(2 - l)} - 1 - \frac{\kappa}{2 - l}.$$

3.3. *Compound Poisson-process driven Ornstein–Uhlenbeck processes.* Let $X$ be an Ornstein–Uhlenbeck process driven by a finite rate subordinator:

$$dX_t = -\mu X_t \, dt + dZ_t,$$

where $Z_t := \sum_{i=1}^{N_t} W_i$, $\{W_i\}_{i \geq 1}$ is an independent and identically distributed collection of random variables from probability measure $F$ and $\{N_t\}_{t \geq 0}$ is a Poisson-process of finite rate $\lambda$, independent of the collection $\{W_i\}_{i \geq 1}$. Such processes are used as storage models (see, e.g., [18]) and have recently been used in financial econometrics as models for stochastic volatility (see [1]).

The exponential decay of $X$, except at jump points, leads to geometric ergodicity of $X$ when the tails of $F(\cdot)$ are sufficiently light. Here we shall explore the case where $F(\cdot)$ is extremely heavy-tailed. First we make this concept precise: we say a probability measure is *heavy-tailed* if, under that probability measure, for all $\kappa > 0$, $\mathbf{E}[e^{\kappa X}] = \infty$. Now let $G$ denote the law of the log jump sizes, that is, $G(A) = F(e^A)$. We have the following negative result showing that for sufficiently heavy-tailed jumps, geometric ergodicity and even ergodicity can fail. As usual, we let $\pi$ denote the invariant probability measure (should it exist).

LEMMA 17. (i) *Suppose $\int x G(dx) = \infty$, then $X$ fails to be positive recurrent.*

(ii) *Suppose $G$ is heavy-tailed, then $X$ fails to be geometrically ergodic.*

PROOF. Suppose $X_0 = 2$ and consider the petite set $C = [0, 1]$. Then

$$\mathbf{P}(\tau_C > t) \geq \mathbf{P}[\text{jump of size} \geq e^{\mu t} \text{ occurs before time } \log 2/\mu]$$

$$\tag{23} = (1 - 2^{-\lambda/\mu}) \int_{\mu t}^{\infty} G(x) \, dx.$$

For positive recurrence, we require that $\mathbf{E}(\tau_C)$ to be finite, that is, that $\mathbf{P}[\tau_C > t]$ be integrable. However, the integrated right-hand side of (23) is just

$$\int_0^\infty dt \int_{\mu t}^\infty G(x)\,dx = \int_0^\infty \mu^{-1} x G(x)\,dx = \infty$$

by hypothesis, so that $\mathbf{E}(\tau_C) = \infty$ too, so that positive recurrence must fail, proving (i). For (ii), we recall that for geometric ergodicity, we require that for some $\kappa > 0$, $\mathbf{E}[e^{\kappa \tau_C}] < \infty$. (Although not necessary, we shall again assume that $X_0 = 2$ and $C = [0,1]$.) Thus, from (23) we require that

$$(24) \qquad \int_0^\infty e^{\kappa t}\,dt \int_{\mu t}^\infty G(x)\,dx = \mu \kappa^{-1} \int_0^\infty (e^{\kappa x} - 1) G(\mu x)\,dx < \infty.$$

However, this is precluded by the heavy-tailed nature of $G$, thus proving (ii). □

Examples of jump distributions for which geometric ergodicity fails (case 2 above), though we will see that $X$ is positive Harris-recurrent, include the following:

$$F(dx) = \frac{dx}{x(\log x)^k} \qquad \text{at least for } k > 1;$$

$$F(dx) = \frac{e^{-(\log x)^\beta}\,dx}{x} \qquad \text{for some } \beta \le 1.$$

LEMMA 18.  *Suppose that for some* $r > 1$, $m_r := \int_0^\infty [\log(1+u)]^r F(du) < \infty$. *Then,* $X$ *is polynomially ergodic with rate* $(1+t)^{(r-1)}$.

PROOF.  For differentiable functions $V$ in the domain of $\mathcal{A}$,

$$\mathcal{A}V = \int_0^\infty (V(x+u) - V(x))\lambda F(du) - \mu x V'(x).$$

Now set $V(x) = (\log x)^r$, then by direct calculation,

$$(25) \quad \mathcal{A}V^\eta = \int_0^\infty ((\log(x+u))^{r\eta} - (\log x)^{r\eta})\lambda F(du) - \frac{\mu x r \eta (\log x)^{\eta r - 1}}{x}.$$

Now the finiteness of $m_r$ merely ensures the finiteness of the first term on the right-hand side of (25). So, noting that $(\log x)^{r\eta}$ is concave beyond $x = e^{r-1}$ for all $0 < \eta \le 1$, we find that, in fact, the first term on the right-hand side of (25) is bounded as a function of $x$, so that for some positive constant $c$,

$$\mathcal{A}V^\eta \le \int_0^\infty \frac{ur(\log x)^{\eta r - 1}}{x}\lambda F(du)c - r\mu(\log x)^{\eta r - 1}.$$

It is easy to check that all bounded sets are petite in this example, and, therefore, the conditions for the application of Corollary 6 with $\alpha = r^{-1}$. □

**4. Proofs of Section 2.** When not explicitly defined, $c$ denotes a generic finite positive constant. $\theta$ is the usual shift operator on the canonical probability space of the strong Markov process.

LEMMA 19. *If $\Psi^{-1}$ is a Young function and $r \in \Lambda_0$ (resp. $\Lambda$), $[\Psi(r) \vee 1] \in \Lambda_0$ (resp. $\Lambda$).*

PROOF. Let $r \in \Lambda_0$. $\Psi^{-1}$ is a continuous, increasing and convex function, so $\Psi$ is measurable and bounded on bounded sets ([16], Chapter 1). Furthermore, there exists a right-continuous nondecreasing function $\phi$ such that $\ln \Psi(r(t)) = \ln r(t) + \ln\{r(t)^{-1} \int_0^{r(t)} \phi(s)\, ds\}$; thus proving that $\ln \Psi(r(t))/t \downarrow 0$ as $t \to \infty$. This yields $\Psi \in \Lambda_0$. The second assertion deduces easily from the definition of $\Lambda$ and the upper bound $\sup_{t \geq 1} \Psi(at)/\Psi(t) < \infty$ for all $a > 0$ ([16], Chapter 1, pages 7 and 8). □

While Theorem 1 and Corollary 6 are claimed for a rate function $r \in \Lambda$, Lemma 19 shows that they can be established for a rate $r \in \Lambda_0$, and we will do so.

4.1. *Proof of Theorem 1.* Without loss of generality, we assume $\Psi_1(r_*) \geq \mathbb{1}$ and $\Psi_2(f_*) \geq 1$.

LEMMA 20. *Let $r \in \Lambda_0$ and $f \geq 1$ be a Borel function. For any closed set $C$ such that $\sup_C G_C(\cdot, f, r; \delta) < \infty$, there exists a constant $M < \infty$ such that for all $x \in \mathcal{X}$ and $t \geq \delta$, $G_C(x, f, r; t) \leq M^{\lfloor t/\delta \rfloor} G_C(x, f, r; \delta)$.*

PROOF. The proof is on the same lines as the proof of Lemma 4.1 in [20] that addresses the case $r = \mathbb{1}$, and the details are omitted. Using the property $r(s+t) \leq r(s)r(t)$ ([34], Lemma 1(d)), we obtain $M = 1 + \sup_{t \geq \delta}[r(t)/r^0(t)] \times \sup_C G_C(\cdot, f, r; \delta)$, which is finite since $\lim_t r(t)/r^0(t) = 0$ (this is a consequence of [34], Lemma 1). □

PROPOSITION 21. *Let $r \in \Lambda_0$ and $f \geq 1$ be a Borel function. Assume that $X$ is $\phi$-irreducible and $\sup_C G_C(\cdot, f, r; \delta) < \infty$ for some closed petite set $C$ and $\delta > 0$. $x \mapsto G_C(x, f, r; \delta)$ is finite $\psi$-almost surely for some (and then any) maximal irreducibility measure $\psi$, and $C$ is accessible.*

PROOF. By [20], Proposition 3.2(ii), for all $\lambda > 0$, there exist a positive integer $m$ and a maximal irreducibility measure $\psi$ such that $\psi(\cdot) \leq \inf_{x \in C} R_\lambda^m(x, \cdot)$, where $R_\lambda$ is the resolvent kernel $R_\lambda(x, \cdot) = \int \lambda \exp(-\lambda t) P^t(x, \cdot)\, dt$. By Lemma 20, $R_\lambda G_C(\cdot, f, r; \delta)(x) \leq c G_C(x, f, r; \delta)$, where $c$ is finite for some convenient $\lambda$. Hence, $\psi G_C(\cdot, f, r; \delta) < \infty$, proving the first statement. This implies that there exists an accessible set $B$ such that $\sup_{x \in B} \mathbb{E}_x[\tau_C(\delta)] \leq \sup_{x \in B} G_C(x, f, r; \delta) < \infty$. Then for $q$ large enough, $\inf_{x \in B} \mathbb{P}_x(\tau_C(\delta) \leq q) > 0$ and, for any $x$, $\mathbb{E}_x[\eta_C] \geq P^n(x, B) \inf_{x \in B} \mathbb{P}_x(\tau_C(\delta) \leq q) > 0$ for some $n$ depending upon $(x, B)$. Hence, $C$ is accessible.  $\square$

PROPOSITION 22.  *Suppose assumptions* (i) *and* (ii) *of Theorem* 1. *Then:*

(i)  *There exist $t_0$ and a measure $\nu$ such that $\inf_{t \geq t_0} \inf_{x \in C} P^t(x, \cdot) \geq \nu(\cdot)$, and $\nu(C) > 0$.*
(ii)  *For any set $B$ such that $\nu(B) > 0$, $\mathbb{E}_x[r^0(T_{m,B})] \leq R_B \mathbb{E}_x[r^0(\tau_C(\delta))]$ for some finite constant $R_{t,B}$.*
(iii)  *For any $t \geq 0$ and any accessible set $B$, $\mathbb{E}_x[r^0(\tau_B(t))] \leq R_{t,B} \mathbb{E}_x[r^0(\tau_C(\delta))]$ for some finite constant $R_{t,B}$.*

PROOF.  (i) Results from R2, Proposition 21 and Lemma 2. (ii) Let $t_0$ and $\nu$ be given by (i). Set $\tau = \tau_C(t_0 + m)$; and define the sequence of iterates $\tau^1 = \tau$ and for $n \geq 2$, $\tau^n = \tau^{n-1} + \tau \circ \theta^{\tau^{n-1}}$. Finally, let $(u_n)_{n \geq 2}$ be a $\{0, 1\}$-valued process given by $u_n = 1$ if $X_{\lceil (\tau^{n-1} + t_0)/m \rceil m} \in B$ and 0 otherwise. $\lceil t \rceil$ denotes the upper integer part of $t$. Then $u_n \in \mathcal{H}_n$ with $\mathcal{H}_n = \sigma(X_t, t \leq \tau^n)$, and by the strong Markov property, $\mathbb{P}_x(u_n = 1 | \mathcal{H}_{n-1}) \geq \nu(C) > 0$ for $n \geq 2$. Finally, set $\eta = \inf\{n \geq 2, u_n = 1\}$, so that $\mathbb{E}_x[r^0(T_{m,B})] \leq \mathbb{E}_x[r^0(\tau^\eta)]$. Using again the strong Markov property and the inequality $r^0(t_1 + t_2) \leq r^0(t_1) + r(t_1)r^0(t_2)$ [34],

$$\mathbb{E}_x[r^0(\tau^\eta)] \leq \sum_{n \geq 2} \mathbb{E}_x[r^0(\tau^n) \mathbb{1}_{\eta \geq n}] = \sum_{n \geq 2} \left\{ a_x(n) + \sup_{x \in C} \mathbb{E}_x[r^0(\tau)] b_x(n) \right\}, \tag{26}$$

for all $n \geq 2$, where $a_x(n) = \mathbb{E}_x[r^0(\tau^{n-1}) \mathbb{1}_{\eta \geq n}]$ and $b_x(n) = \mathbb{E}_x[r(\tau^{n-1}) \mathbb{1}_{\eta \geq n}]$. Since, by Lemma 20, $\sup_C \mathbb{E}_x[r^0(\tau)] < \infty$, there exists $0 < \rho < 1$ and a finite constant $c$ such that

$$b_x(n) \leq \rho b_x(n-1) + c(1 - \nu(C))^{n-1},$$
$$a_x(n) \leq (1 - \nu(C)) a_x(n-1) + b_x(n-1) \sup_{x \in C} \mathbb{E}_x[r^0(\tau)];$$

and $b_x(2) = \mathbb{E}_x[r(\tau)]$, $a_x(2) = \mathbb{E}_x[r^0(\tau)]$. The proof is on the same lines as the proof of [25], Lemma 3.1, and is omitted for brevity. Hence, $\mathbb{E}_x[r^0(T_{m,B})] \leq c(\mathbb{E}_x[r^0(\tau)] + \mathbb{E}_x[r(\tau)])$ for some $c < \infty$. The proof is concluded, applying again Lemma 20 and the bound $\sup_{t \geq a} r(t)/r^0(t) < \infty$ for all $a > 0$ (see the proof of Lemma 19).

(iii) $B$ is accessible and $C$ petite so there exist $t_0 \geq 0$ and $\gamma > 0$ such that $\inf_{x \in C} \mathbb{P}_x(\tau_B \leq t_0 + t) \geq \inf_{x \in C} \mathbb{P}_x(\tau_B \leq t_0) \geq \gamma$. Set $\tau = \tau_C(t + t_0)$ and $u_n = 1$ if for some $\tau^{n-1} \leq s \leq \tau^{n-1} + t + t_0$, $X_s \in B$; and $u_n = 0$ otherwise. Following the same lines as in the proof of (ii), it may be proved that there exists $c < \infty$ such that $\mathbb{E}_x[r^0(\tau_B(t))] \le c\mathbb{E}_x[r^0(\tau_C(t + t_0))]$. The proof is concluded by applying Lemma 20. □

PROPOSITION 23. *Suppose assumptions* (i) *and* (ii) *of Theorem* 1. *For any* $(\Psi_1, \Psi_2) \in \mathcal{I}$, $C$ *is a* $(\Psi_2(f_*), \Psi_1(r_*))$-regular set for the process, that is, $\sup_C G_B(\cdot, \Psi_2(f_*), \Psi_1(r_*); t) < \infty$ for any $t > 0$ and any accessible set $B$. $G_B(x, \Psi_2(f_*), \Psi_1(r_*); t) < \infty$ for all $x \in \mathcal{S}_\Psi$ and $\pi(\mathcal{S}_\Psi) = 1$.

PROOF. $(\Psi_2(f_*), \Psi_1(r_*))$-regularity is a consequence of Young's inequality (4), the $(f_*, \mathbb{1})$-regularity of $C$ ([20], Proposition 4.1) and Proposition 22(iii). For the second statement, write

$$G_B(x, \Psi_2(f_*), \Psi_1(r_*); t)$$

$$\le G_C(x, \Psi_2(f_*), \Psi_1(r_*); t) + \mathbb{E}_x\left[\int_{\tau_C(t)}^{\tau_B \circ \theta^{\tau_C(t)}} \Psi_1(r_*(s))\Psi_2(f_*(X_s)) \, ds\right].$$

The result now follows from the strong Markov property, Lemma 20 and the inequality $\Psi_1(r_*(s+t)) \le \Psi_1(r_*(s))\Psi_1(r_*(t))$, which holds since $\Psi_1 \circ r_* \in \Lambda_0$. Finally, $\pi(\mathcal{S}_\Psi) = 1$ by Proposition 21. □

PROPOSITION 24. *Suppose assumptions* (i) *and* (ii) *of Theorem* 1. *The skeleton chain* $P^m$ *is* $\psi$-*irreducible and aperiodic and possesses an accessible petite set* $A$ *such that for all* $(\Psi_1, \Psi_2) \in \mathcal{I}$,

$$(27) \qquad \sup_{x \in A} \mathbb{E}_x\left[\sum_{k=0}^{T_{m,A}-1} \Psi_1(r_*(k))\Psi_2(f_*(X_{km}))\right] < \infty.$$

PROOF. For the definitions of accessibility, smallness, petiteness and aperiodicity of a discrete-time Markov chains, see [21]. From Proposition 22(i), $C$ is small for the skeleton $P^m$ and the skeleton is aperiodic (Theorem 5.4.4, [21]). In addition, by R2, the skeleton is positive and $\pi(f_*) < \infty$. Let $C_n$ be a petite set (for the skeleton $P^m$) such that $A = C \cap C_n$ is of positive $\nu$-measure

$$(28) \qquad \sup_{x \in C_n} \mathbb{E}_x\left[\sum_{k=0}^{T_{m,B}-1} f_*(X_{km})\right] < \infty,$$

for any accessible set $B$ (for the skeleton); the existence of such a set is a consequence of Theorems 14.2.3 and 14.2.11 in [21] and Proposition 22(ii). The set $A$ is accessible and petite for the skeleton. (27) now results from Young's inequality (4), (28) and Proposition 22(ii). □

PROOF OF THEOREM 1. By Proposition 24 and [35], Theorem 2.1 and Proposition 3.2, $\lim_{n\to\infty} \Psi_1(r_*(n)) \| P^{nm}(x,\cdot) - \pi(\cdot) \|_{\Psi_2(f_*)} = 0$ for $\pi$ a.a. $x$. By Jensen's inequality, the upper bound $\sup_{t\geq 1} \Psi_2(at)/\Psi_2(t) < \infty$ for all $a > 0$, and assumption (iii), we have for all $t \leq m$, $P^t \Psi_2(f_*) \leq c\Psi_2(f_*)$. In addition, since $\Psi_1(r_*) \in \Lambda$, $\Psi_1(r_*(n+t)) \leq c\Psi_1(r_*(n))$ for all $t \leq m$ ([34], Lemma 1). Hence,

$$\tag{29} \lim_{t\to\infty} \Psi_1(r_*(t)) \| P^t(x,\cdot) - \pi(\cdot) \|_{\Psi_2(f_*)} = 0, \qquad \pi \text{ a.a. } x.$$

We now prove that this convergence occurs for all $x \in \mathcal{S}_\Psi$ which is of $\pi$-measure one, by Proposition 23. To that goal, we mimic the proof of [22], Theorem 7.2. By Egorov's theorem, there exists a set $A$, $\pi(A) > 0$, such that (29) holds uniformly for all $x \in A$. For all Borel functions, $g \in \mathcal{L}_{\Psi_2(f_*)}$, set $\bar{g} := g - \pi(g)$. Since $\Psi_1(r_*) \in \Lambda_0$,

$$\Psi_1(r_*(t)) |\mathbb{E}_x[\bar{g}(X_t)\mathbb{1}_{\tau_A \leq t}]| \leq \Psi_1(r_*(t)) \int_0^t \sup_{y \in A} |P^{t-s}\bar{g}|(y) \mathbb{P}_x(\tau_A \in ds)$$
$$\leq M\{\mathbb{E}_x[\Psi_1(r_*(\tau_A))] + \mathbb{E}_x[\Psi_1(r_*(\tau_A))\mathbb{1}_{\tau_A \geq t/2}]\},$$

where $M = \sup_{y \in A} \sup_{s \geq 0} r_*(s) |P^s \bar{g}|(y)$. Let $x \in \mathcal{S}_\Psi$; from Proposition 23, $\mathbb{E}_x[\Psi_1(r_*(\tau_A))] < \infty$ and $\lim_{t\to\infty} \mathbb{E}_x[\Psi_1(r_*(\tau_A))\mathbb{1}_{\tau_A \geq t/2}] = 0$. Since the limit (29) holds uniformly for all $x \in A$, $M$ is finite. Hence, $\lim_{t\to\infty} \Psi_1(r_*(t)) |\mathbb{E}_x[\bar{g}(X_t) \times \mathbb{1}_{\tau_A \leq t}]| = 0$ uniformly for all $g \in \mathcal{L}_{\Psi_2(f_*)}$.

Since $\pi(f_*) < \infty$, $|\mathbb{E}_x[\bar{g}(X_t)\mathbb{1}_{\tau_A \geq t}]| \leq c\mathbb{E}_x[\Psi_2(f_*(X_t))\mathbb{1}_{\tau_A \geq t}]$. Following the same lines as in the proof of [22], Theorem 7.2, using again $\sup_{u \leq m} P^u f_* \leq cf_*$, we obtain

$$\Psi_1(r_*(t)) \mathbb{E}_x[\Psi_2(f_*(X_t))\mathbb{1}_{\tau_A \geq t}]$$
$$\leq c\Psi_1(r_*(m)) \inf_{0 \leq u \leq m} \Psi_1(r_*(t-u)) \mathbb{E}_x[\Psi_2(f_*(X_{t-u}))\mathbb{1}_{\tau_A \geq t-u}].$$

By Proposition 23, $G_A(x, \Psi_2(f_*), \Psi_1(r_*); 0) < \infty$, which implies that the upper limit in the right-hand side is zero, proving that $\lim_{t\to\infty} \Psi_1(r_*(t)) \times \mathbb{E}_x[\Psi_2(f_*(X_t))\mathbb{1}_{\tau_A \geq t}] = 0$. Hence, uniformly for $g \in \mathcal{L}_{\Psi_2(f_*)}$, $\lim_{t\to\infty} \Psi_1(r_*(t)) \times |\mathbb{E}_x[\bar{g}(X_t)\mathbb{1}_{\tau_A \geq t}]| = 0$. This concludes the proof. $\square$

4.2. *Proof of Corollary* 6. Set $f_* := V^{1-\alpha}$ and $r_*(t) := (t+1)^{1/\alpha-1}$.

LEMMA 25. *Suppose assumption* (ii) *of Corollary* 6. *For any* $\alpha \leq \eta \leq 1$, $t \geq 0$, *and any* $\mathcal{F}_t$-*stopping-time* $\tau$,

$$c_\eta \mathbb{E}_x\left[\int_0^{\tau \wedge t} V^{\eta-\alpha}(X_s)\,ds\right] + \mathbb{E}_x[V^\eta(X_{\tau \wedge t})] \leq V^\eta(x) + b\mathbb{E}_x\left[\int_0^{\tau \wedge t} \mathbb{1}_C(X_s)\,ds\right].$$

PROOF.  By definition of $\mathcal{A}V$,

$$c_\eta \mathbb{E}_x \left[ \int_0^{\tau \wedge t \wedge T_n} V^{\eta-\alpha}(X_s) \, ds \right] + \mathbb{E}_x[V^\eta(X_{\tau \wedge t \wedge T_n})]$$

$$\leq V^\eta(x) + b\mathbb{E}_x \left[ \int_0^{\tau \wedge t \wedge T_n} \mathbb{1}_C(X_s) \, ds \right].$$

The right-hand side is upper bounded by $V(x) + bt$ and by the monotone convergence theorem, it converges to $V(x) + b\mathbb{E}_x[\int_0^{\tau \wedge t} \mathbb{1}_C(X_s) \, ds]$ as $n \to \infty$. The lemma now results from Fatou's lemma. □

PROPOSITION 26.  *Suppose assumption* (ii) *of Corollary* 6. *For all* $\delta > 0$, *there exists* $c < \infty$ *such that for all* $x \in \mathcal{X}$, $G_C(x, \mathbb{1}, r_*; \delta) \leq cV(x)$.

PROOF.  Set $q := \lfloor 1/\alpha \rfloor$, where $\lfloor \cdot \rfloor$ denotes the lower integer part. By Lemma 25, we have $\mathbb{E}_x[\tau_C] \leq cV^\alpha(x)$ and by Jensen's inequality, we obtain $\mathbb{E}_x[\tau_C^{\alpha^{-1}-q}] \leq cV^{1-q\alpha}(x)$. We prove by induction that for all integer $1 \leq l \leq q$, $\mathbb{E}_x[\tau_C^{\alpha^{-1}-l}] \leq cV^{1-l\alpha}(x)$. The case $l = q$ holds; assume it is verified for some $2 \leq l \leq q$. The induction hypothesis and Lemma 25 yield

$$\mathbb{E}_x[\tau_C^{\alpha^{-1}-l+1}] \leq c\mathbb{E}_x \left[ \int_0^{\tau_C} \mathbb{E}_{X_s}[\tau_C^{\alpha^{-1}-l}] \, ds \right]$$

$$\leq c\mathbb{E}_x \left[ \int_0^{\tau_C} V^{1-l\alpha}(X_s) \, ds \right] \leq cV^{1-l\alpha+\alpha}(x),$$

which concludes the induction. For $l = 1$, this yields $G_C(x, \mathbb{1}, r_*; 0) \leq cV(x)$. Finally, by standard manipulations and Lemma 25, we have $G_C(x, \mathbb{1}, r_*; \delta) \leq c(1 + P^\delta V(x)) \leq cV(x)$. □

PROOF OF COROLLARY 6.  We check the conditions for the application of Theorem 1. Lemma 25 and Proposition 26 imply $G_C(x, f_*, ; \delta) \leq cV(x)$ and $G_C(x, \mathbb{1}, r_*; \delta) \leq cV(x)$, from which we deduce the condition (ii) of Theorem 1, and by R1, condition (i) of Theorem 1. Condition (iii) follows from Lemma 25. Finally, $\mathcal{S}_\Psi = \mathcal{X}$. □

**Acknowledgment.**  We are grateful to the referee for his comments, and, in particular, for bringing [2] to our attention.

CNRS/LMC-IMAG
51 rue des Mathematiques
BP 53
38041 Grenoble Cedex 9
France
e-mail: Gersende.Fort@imag.fr

Department of Mathematics
  and Statistics
Lancaster University
Lancaster LA1 4YF
United Kingdom
e-mail: g.o.roberts@lancaster.ac.uk